\documentclass[journal ]{new-aiaa}
\usepackage[utf8]{inputenc}
\usepackage{textcomp}
\usepackage{ulem} 
\usepackage{graphicx}
\usepackage{amsmath}
\usepackage[version=4]{mhchem}
\usepackage{siunitx}
\usepackage{longtable,tabularx}
\usepackage{color}
\usepackage{mathtools}
\usepackage{float}
\usepackage{caption,subcaption}
\usepackage{tikz}
\usepackage[many]{tcolorbox}
\usepackage{color}
\usepackage{epstopdf}
\usepackage[numbers,sort&compress]{natbib}
\usepackage{algorithm}
\usepackage{algorithmic}
\usepackage{setspace}

\def\1{\boldsymbol{1}}

\title{A New Smoothing Technique for Bang-Bang Optimal Control Problems}

\author{Kun Wang\footnote{PhD student, School of Aeronautics and Astronautics. email: \underline{wongquinn@zju.edu.cn}.}, Zheng Chen\footnote{Researcher, School of Aeronautics and Astronautics. email: \underline{z-chen@zju.edu.cn} (Corresponding author). }, Zhenyu Wei\footnote{PhD student, College of Control Science and Engineering. }, Fangmin Lu\footnote{PhD student, School of Aeronautics and Astronautics.}, and Jun Li\footnote{Professor, School of Aeronautics and Astronautics.}}

\affil{ Zhejiang University, Hangzhou 310027, Zhejiang, China}

\begin{document}

\maketitle

\begin{abstract}
Bang-bang control is ubiquitous for Optimal Control Problems (OCPs) where the constrained control variable appears linearly in the dynamics and cost function. Based on the Pontryagin’s Minimum Principle, the indirect method is widely used to numerically solve OCPs because it enables one to derive the theoretical structure of the optimal control. However, discontinuities in the bang-bang control structure may result in numerical difficulties for the gradient-based indirect method. In this case, smoothing or regularization procedures are usually applied to eliminating the discontinuities of bang-bang controls. Traditional smoothing or regularization procedures generally modify the cost function by adding a term depending on a small parameter, or introducing a small error into the state equation. Those procedures may complexify the numerical algorithms or degenerate the convergence performance. To overcome these issues, we propose a bounded smooth function, called normalized $L^2$-norm function, to approximate the sign function in terms of the switching function. The resulting optimal control is smooth and can be readily embedded into the indirect method. Then, the simplicity and improved performance of the proposed method over some existing methods are numerically demonstrated by a minimal-time oscillator problem and a minimal-fuel low-thrust trajectory optimization problem that involves many revolutions. 
\end{abstract}

\section{ Introduction}
\label{intro}
Optimal control theory is widely applied to determining a sequence of inputs that steer a dynamical system from an initial condition to a desired final condition while minimizing a specified cost function. If the dynamical system is nonlinear so that the analytical solution cannot be obtained, numerical algorithms are often needed to find the optimal solution. Generally speaking, an Optimal Control Problem (OCP) can be numerically solved via either the direct or indirect method as summarized in Ref. \cite{rao2009survey}. The indirect method transforms an OCP into a Two-Point Boundary-Value Problem (TPBVP) by applying the necessary conditions for optimality derived from the Pontryagin's Minimum Principle (PMP), which can be leveraged to reveal the structure of optimal control \cite{Bryson:75}.

Bang-bang control frequently arises when the dynamics and cost function are both linear w.r.t. the control variable constrained by an upper and lower bound. Conventionally, there are two types of cost function leading to bang-bang controls. The first one involves $L^0$-norm in terms of the control variable, i.e., a minimal time control problem; see, e.g., Refs. \cite{mall2017epsilon,wijayatunga2023exploiting}. The second type consists of $L^1$-norm dependent on the control variable. In this case, the allowable control is generally nonnegative, resulting in a minimal fuel problem; see, e.g., Refs. \cite{chen20161,wijayatunga2023exploiting}. If gradient-based methods, such as Newton- or quasi-Newton methods, are used to solve the resulting TPBVP, it is quite challenging because of the discontinuities arising from bang-bang controls \cite{bertrand2002new}.

To circumvent the aforementioned numerical difficulties, multiple shooting method was developed in the literature by separating the bang-bang control trajectory into different arcs, resulting in a Multiple-Point Boundary-Value Problem (MPBVP) \cite{mall2020uniform}. To this end, precise structure of the optimal control profile is usually needed {\it a prior} to make the MPBVP tractable, which is almost impossible for highly nonlinear OCPs. Another popular approach is to regularize the OCP by introducing a perturbed energy term into the cost function \cite{bertrand2002new,silva2010smooth} or the state equation \cite{mall2017epsilon}. Specifically, in the work \cite{bertrand2002new} by Bertrand and Epenoy, a barrier function was added to the original cost function, leading to a continuously differentiable optimal control solution. On the other hand, by introducing the concept of error controls and an error parameter, a smoothing procedure for minimal time control problem was developed in Ref. \cite{silva2010smooth}. In Ref. \cite{mall2017epsilon}, a new method was established using the concept of trigonomerization, and the error parameter was introduced only into one state equation, making it easier to implement than the procedure in Ref. \cite{silva2010smooth}.

Recently, a smoothing technique based on the hyperbolic tangent function was proposed to smooth bang-bang controls in Ref. \cite{taheri2018generic}. The implementation of that technique is to simply replace the discontinuously bang-bang control profile with the hyperbolic tangent function dependent on the switching function and a smoothing parameter. Consequently, there is no need to modify the cost function or state equation for such technique. Additionally, it was demonstrated by a minimum-fuel low-thrust orbital transfer problem in Ref. \cite{taheri2018generic} that using the hyperbolic tangent function to approximate bang-bang controls can significantly reduce the computational time.  

In this paper, we propose a new smoothing technique, which uses a bounded smooth function to approximate the sign function in terms of the switching function. The bounded smooth function, called normalized $L^2$-norm function, is expressed as $\frac{S}{\sqrt{\delta+\lvert S\rvert^2}}$, where $S$ is the switching function derived from the PMP, and $\delta$ is a smoothing constant. Similar to the work in Ref. \cite{taheri2018generic}, after the switching function is derived,
our method can be implemented by simply replacing the bang-bang control profile with the normalized $L^2$-norm function. 

The structure of this paper is as follows. Sect.~\ref{Pre} reviews the necessary conditions for optimality for a nonlinear control-affine system where the optimal control takes a bang-bang form. The smoothing technique is presented in Sect.~\ref{Approximating}. In Sect.~\ref{Test},
numerical simulations on two test problems, including 
a minimal-time oscillator and a minimal-fuel low-thrust transfer involving many revolutions, are used to illustrate the benefits of the proposed method. 
\section{Preliminaries}
\label{Pre}
\subsection{Problem Statement}
We consider a nonlinear control-affine system as
\begin{align}
\dot{\boldsymbol{x}}(t) = \boldsymbol{f}(\boldsymbol{x},\boldsymbol{u},t) =
\boldsymbol{f_0}(\boldsymbol{x},t) + \boldsymbol{f_1}(\boldsymbol{x},t)\boldsymbol{u},
\label{EQ:dyna_equation}
\end{align}
where $\boldsymbol{x} \in \mathbb{R}^n$, $\boldsymbol{u} \in \mathbb{R}^m$, and $t \in \mathbb{R}_0^+$ represent the state vector, control vector, and time, respectively; $\boldsymbol{f} :\mathbb{R}^n \times \mathbb{R}^m \times \mathbb{R}_0^+ \rightarrow \mathbb{R}^n$, 
$\boldsymbol{f_0} :\mathbb{R}^n \times \mathbb{R}_0^+ \rightarrow \mathbb{R}^n$, and $\boldsymbol{f_1} :\mathbb{R}^n \times \mathbb{R}_0^+ \rightarrow \mathbb{R}^{n\times m}$ are smooth vector fields.

Without loss of generality, consider to minimize
\begin{align}
J = \int_{t_0}^{t_f}{L(\boldsymbol{x},\boldsymbol{u},t)}dt,
\label{EQ:cost_equation}
\end{align}
where $t_0$ and $t_f$ denote the initial time and final time, respectively; $t_0$ is usually known and $t_f$ can be either free or fixed. The mapping $L:\mathbb{R}^n \times \mathbb{R}^m \times \mathbb{R}_0^+ \rightarrow \mathbb{R}$ is the running cost and is linear w.r.t. the control vector $\boldsymbol{u}$. Note that the cost function $J$ in Eq.~(\ref{EQ:cost_equation}) can be transformed into other forms, such as Mayer and Bolza type \cite{longuski2014optimal}.
The initial constraint is given by
\begin{align}
\boldsymbol{\Psi}(\boldsymbol{x}(t_0),t_0) = \boldsymbol{0},
\label{EQ:initial_equation}
\end{align}
and the terminal constraint is given by
\begin{align}
\boldsymbol{\Phi}(\boldsymbol{x}(t_f),t_f) = \boldsymbol{0},
\label{EQ:terminal_equation}
\end{align}
where $\boldsymbol{\Psi}:\mathbb{R}^n \times \mathbb{R}_0^+ \rightarrow \mathbb{R}^p$, and $\boldsymbol{\Phi}:\mathbb{R}^n \times \mathbb{R}_0^+ \rightarrow \mathbb{R}^q$ are smooth w.r.t. their variables. The control vector $\boldsymbol{u}$ takes values in the set $\boldsymbol{\mathscr{U}}$
\begin{align}
\boldsymbol{\mathscr{U}} = \left\{\boldsymbol{u} | u_i \in \left\{u_i^{min}, u_i^{max}\right\}, i=1,\ldots,m\right\}
\label{EQ:control_equation}
\end{align}
where $u_i^{min}$ and $u_i^{max}$ denote the lower bound and the upper bound of $u_i$, respectively.

Thus, the corresponding OCP is to find a measurable control $\boldsymbol{u}(t) \in \boldsymbol{\mathscr{U}}$ on $[0,t_f]$ that drives the dynamic system described by Eq.~(\ref{EQ:dyna_equation}) from the initial condition in Eq.~(\ref{EQ:initial_equation}) to the terminal condition in Eq.~(\ref{EQ:terminal_equation}) while the cost function in Eq.~(\ref{EQ:cost_equation}) is minimized.
\subsection{Necessary Conditions for Optimality}
Denote by $\boldsymbol{\lambda} \in \mathbb{R}^n$ the costate vector associated with the state vector $\boldsymbol{x}$. Then,  the Hamiltonian $\mathscr{H}: \mathbb{R}^n \times \mathbb{R}^m \times \mathbb{R}^n \times \mathbb{R}_0^+ \rightarrow \mathbb{R} $ for the OCP is defined as \cite{Pontryagin}
\begin{align}
\mathscr{H}(\boldsymbol{x},\boldsymbol{u},\boldsymbol{\lambda},t) = L(\boldsymbol{x},\boldsymbol{u},t) + \boldsymbol{\lambda}^T [\boldsymbol{f_0}(\boldsymbol{x},t) + \boldsymbol{f_1}(\boldsymbol{x},t)\boldsymbol{u}].
\label{EQ:Ham_equation}
\end{align}
The costate is governed by
\begin{align}
\dot{\boldsymbol{\lambda}}(t) = -\frac{\partial{\mathscr{H}}(t)}{\partial{\boldsymbol{x}(t)}}.
\label{EQ:costate_equation}
\end{align}
The transversality condition implies
\begin{align}
\boldsymbol{\lambda}(t_f) = \boldsymbol{v_f}^T\frac{\partial \boldsymbol{\Phi}}{\partial{\boldsymbol{x}}(t_f)},
\label{EQ:costate_final}
\end{align}
where $\boldsymbol{v_f} \in \mathbb{R}^q$ is the constant vector of Lagrangian multipliers, used to adjoin the terminal condition. 
The PMP states that the optimal control will minimize the Hamiltonian $\mathscr{H}$ for all admissible values in $\boldsymbol{\mathscr{U}}$ along an optimal trajectory, i.e., 
\begin{align}
\mathscr{H}(\boldsymbol{x}^*,\boldsymbol{u}^*,\boldsymbol{\lambda}^*,t) \leq \mathscr{H}(\boldsymbol{x}^*,\boldsymbol{u},\boldsymbol{\lambda}^*,t).
\label{EQ:control_law}
\end{align}
Recall that the Hamiltonian $\mathscr{H}$ is linear in the control vector $\boldsymbol{u}$, thus it can be rewritten as 
\begin{align}
\mathscr{H}(\boldsymbol{x},\boldsymbol{u},\boldsymbol{\lambda},t) = \mathscr{H}_0(\boldsymbol{x},\boldsymbol{\lambda},t) + \boldsymbol{S}^T(\boldsymbol{x},\boldsymbol{\lambda},t)\boldsymbol{u},
\label{EQ:Ham_new}
\end{align}
where $\mathscr{H}_0:\mathbb{R}^n \times \mathbb{R}^n \times \mathbb{R}_0^+ \rightarrow \mathbb{R}$ is part of $\mathscr{H}$ independent of $\boldsymbol{u}$, and $\boldsymbol{S}:\mathbb{R}^n \times \mathbb{R}^n \times \mathbb{R}_0^+ \rightarrow \mathbb{R}^m$ denotes the switching functions corresponding to $\boldsymbol{u}$.
Denote by $S_{i} = \frac{\partial{\mathscr{H}}}{\partial{u_i}} (i = 1,\ldots,m)$ the switching function associated with the control $u_i$. Assume that the optimal control does not contain any singular arc, then the optimal control can be derived as
\begin{align}
u_i^*(t) = 
\begin{cases}
u_i^{max}, \text{if}~{S_{i}(t)} < {0}, \\
u_i^{min}, \text{if}~{S_{i}(t)} > {0}.
\label{EQ:u_i}
\end{cases}
\end{align}
In addition, because the Hamiltonian $\mathscr{H}$ does not contain the time $t$ explicitly and if the final time $t_f$ is free, it holds 
\begin{align}
\mathscr{H}(\boldsymbol{x},\boldsymbol{u},\boldsymbol{\lambda},t) \equiv 0.
\label{EQ:optimal_new}
\end{align}
\subsection{Formulation of the TPBVP}
Eqs.~(\ref{EQ:dyna_equation}), (\ref{EQ:costate_equation}), and (\ref{EQ:u_i}) constitute a set of Ordinary Differential Equations (ODEs)  in terms of $\boldsymbol{x}$ and $\boldsymbol{\lambda}$. The OCP is formulated as a TPBVP of the form
\begin{align}
\boldsymbol{\zeta}(\boldsymbol{\lambda}(t_0)) = [\boldsymbol {\Phi}(t_f),\boldsymbol{\lambda}(t_f)-\boldsymbol{v_f}^T\frac{\partial \boldsymbol{\Phi}}{\partial{\boldsymbol{x}}(t_f)}] = \boldsymbol{0},
\label{EQ:shooting_function}
\end{align}
or (in the case of free $t_f$)
\begin{align}
\boldsymbol{\zeta}(\boldsymbol{\lambda}(t_0),t_f) = [\boldsymbol {\Phi}(t_f),\boldsymbol{\lambda}(t_f)-\boldsymbol{v_f}^T\frac{\partial \boldsymbol{\Phi}}{\partial{\boldsymbol{x}}(t_f)},\mathscr{H}(t_f)] = \boldsymbol{0},
\label{EQ:shooting_function}
\end{align}
where $\boldsymbol{\zeta}$, and $\boldsymbol{\lambda}(t_0)$ ($[\boldsymbol{\lambda}(t_0),t_f)]$ if $t_f$ is free) is called the shooting function and shooting variable, respectively. To solve the TPBVP, one should provide an initial guess for the shooting variable. Since the optimal control in Eq.~(\ref{EQ:u_i}) is either at the maximum or minimum, it leads to non-smoothness and discontinuity of the ODEs. As a result, singularity of the Jacobian of the shooting function $\boldsymbol{\zeta}$ frequently arises on a large domain \cite{bertrand2002new}, making it daunting to solve the shooting function using gradient-based numerical algorithms. In the next section, we shall propose a bounded smooth function to smooth the bang-bang control profile.
\section{Smoothing the Bang-Bang Control with Normalized $L^2$-Norm Function}
\label{Approximating}
\begin{figure}[!htp]
\centering
\includegraphics[width = 0.75\linewidth]{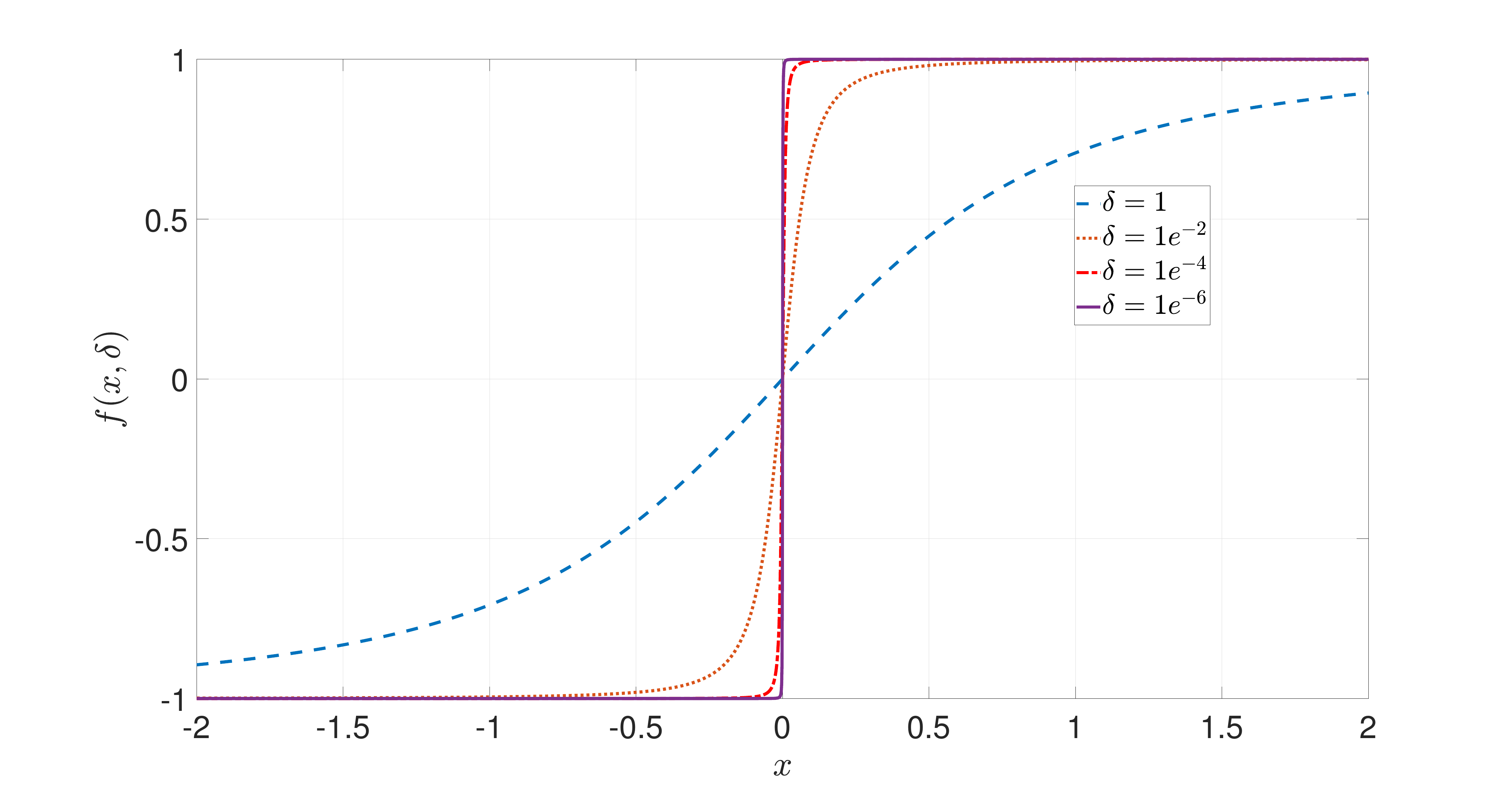}
\caption{Plot of $f(x,\delta)$ versus $x$ for different smoothing constants.}
\label{Fig:norm_function}
\end{figure}

Without loss of generality, assume that the allowable control is a scalar $u \in [u^{min}, u^{max}]$. Then, the optimal control in Eq.~(\ref{EQ:u_i}) can be rewritten as
\begin{align}
u^*(t) =  \frac{1}{2}[(u^{max}+u^{min}) - (u^{max}-u^{min})sgn(S(t))],
\label{EQ:sign_equation}
\end{align}
where $sgn$ is the sign function. In this case, discontinuities of bang-bang controls arise, leading to numerical difficulties for the gradient-based indirect method. To this end, we define a bounded smooth function, called normalized $L^2$-norm function to approximate the sign function, i.e., 
\begin{align}
f(x,\delta):= \frac{x}{\sqrt{\delta+\lvert x\rvert^2}},
\label{EQ:norm_equation_f}
\end{align}
where $x \in \mathbb{R}$ is the input, and $\delta$ is a small positive number, denoting the smoothing constant. Fig.~\ref{Fig:norm_function} shows the normalized $L^2$-norm function over an input $x \in [-2,2]$. 
It can be seen that, $f(x,\delta)$ is bounded by $[-1,1]$, and smaller $\delta$ leads to closer approximation of the sign function. Following the work in Ref. \cite{taheri2018generic}, set the switching function $S(t)$ as the input for $f(x,\delta)$. Then, this function with $\delta \rightarrow 0$ can be employed to approximate the bang-bang control in Eq.~(\ref{EQ:sign_equation}), i.e.,
\begin{align}
u^*(t) \approx  u^*(t,\delta) = \frac{1}{2}[(u^{max}+u^{min}) - (u^{max}-u^{min})\frac{S(t)}{\sqrt{\delta+\lvert S(t)\rvert^2}}].
\label{EQ:norm_equation_approx}
\end{align}

The function defined in Eq.~(\ref{EQ:norm_equation_f}) acts like a filter function that smooths the bang-bang control. In this case, the approximation in Eq.~(\ref{EQ:norm_equation_approx}) can be leveraged to replace the  discontinuous optimal control profile in Eq.~(\ref{EQ:sign_equation}) for the indirect method. In the following section, two test problems will be used to demonstrate the benefits of the proposed method. 
\section{Numerical Simulations}
\label{Test}
In this section, we apply the developed technique to a minimal-time oscillator problem and a minimal-fuel low-thrust trajectory optimization problem. The procedure for establishing the resulting TPBVPs will be briefly discussed first. Then, the subroutine \it{fsolve} (\rm{a Newton-like iterative solver}) will be employed to solve the TPBVPs. A desktop machine equipped with
an Intel(R) Core(TM) i9-10980XE CPU @ 3.0 GHz is used to perform all computations. 
\subsection{Numerical Simulations on a Minimal-Time Oscillator Problem}
\subsubsection{Problem Statement}
Consider an undamped harmonic oscillator problem in Refs. \cite{mall2017epsilon,longuski2014optimal}. The dynamics of the oscillator is 
\begin{align}
\begin{cases}
\dot{x_1}(t) = x_2(t),\\
\dot{x_2}(t) = -x_1(t) + u(t), 
\end{cases}
\label{EQ:dynamics_oscillator}
\end{align}
where $t$ denotes the time, $\boldsymbol{x}(t) = [x_1(t), x_2(t)]^T$ the state vector, and $u(t)$ the input bounded by $\lvert u(t) \rvert \leq 1~\rm{for}~a.e.~ \it{t \in }~[\rm{0},\it{t_f}]$. The initial condition at $t = 0$ is given by 
\begin{align}
x_1(0) = x_{1_0}, x_2(0) = x_{2_0},
\label{EQ:initial_oscillator}
\end{align}
and the final condition is given by
\begin{align}
x_1(t_f) = x_{1_{f}}, x_2(t_f) = x_{2_{f}}.
\label{EQ:final_oscillator}
\end{align}
The cost function to be minimized is defined as
\begin{align}
J = \int_{0}^{t_f} 1~ dt.
\end{align}

The Hamiltonian can be formulated as \cite{Pontryagin} 
\begin{align}
\mathscr{H}(t) = \lambda_1(t)x_2(t) + \lambda_2(t)[-x_1(t) + u(t)] + 1,
\label{EQ:Ham1}
\end{align}
where $\boldsymbol{\lambda}(t) = [\lambda_1(t), \lambda_2(t)]^T$ is the costate vector associated with the state vector $\boldsymbol{x}(t)$. Then, the costate is governed by
\begin{align}
\begin{cases}
\dot{\lambda_1}(t) = -\frac{\partial \mathscr{H}(t)}{\partial x_1(t)} = \lambda_2(t),\\
\dot{\lambda_2}(t) = -\frac{\partial \mathscr{H}(t)}{\partial x_2(t)} = -\lambda_1(t).
\label{EQ:costate_P1}
\end{cases}
\end{align}
Since $\mathscr{H}(t)$ is a linear function of $u(t)$, the switching function can be readily obtained as
\begin{align} 
S(t)=\frac{\partial \mathscr{H}(t)}{\partial u(t)} = \lambda_2(t).
\label{EQ:switching_P1}
\end{align}
Note that a singular arc was proved non-optimal in Ref. \cite{longuski2014optimal}, thus the optimal control that minimizes $\mathscr{H}(t)$ is
\begin{align}
u^*(t) = 
\begin{cases}
1,\text{if}~{\lambda}_2(t) < \rm{0}, \\
-1,\text{if}~{\lambda}_2(t) > \rm{0}.
\label{EQ:u_P1}
\end{cases}
\end{align}
Since the final time $t_f$ is free, along an optimal solution, it holds
\begin{align}
\mathscr{H}(t) \equiv 0.
\label{EQ:HamFinal_P1}
\end{align}

Thanks to the normalized $L^2$-norm function, the discontinuous optimal control in Eq.~(\ref{EQ:u_P1}) can be approximated by
\begin{align}
u^*(t) \approx u^*(t,\delta) =  -\frac{S(t)}{\sqrt{\delta+\lvert S(t)\rvert^2}}.
\label{EQ:smoothing_P1}
\end{align}
The resulting TPBVP amounts to solving the shooting function
\begin{align}
\boldsymbol{\zeta}(\lambda_1(0), \lambda_2(0),t_f) = [x_1{(t_f)}-x_{1_f},x_2{(t_f)}-x_{2_f}, \mathscr{H}(t_f)] = \boldsymbol{0},
\label{EQ:shooting_functionp1}
\end{align}
where $[\lambda_1(0), \lambda_2(0),t_f]^T$ 
is the shooting variable.
\subsubsection{Numerical Results}
We use the same boundary conditions specified in Ref. \cite{mall2017epsilon}, and they are presented in Table~\ref{Table:Boundarys}.
\begin{table}[!htp]
\centering
\caption{Boundary conditions for the minimal-time oscillator problem}
\begin{tabular}{ccccc}
\hline
Parameter  &$x_{1_0}$ &$x_{2_0}$  &$x_{1_f}$  &$x_{2_f}$ \\
\hline
Value     &$1$    &$1$  &$0$ &$0$ \\
\hline
\label{Table:Boundarys}
\end{tabular}
\end{table}
To demonstrate the performance of the proposed method, following the work in Ref. \cite{taheri2018generic} by Taheri and Junkins, $10,000$ simulations are performed for three smoothing methods, including the epsilon-trig regularization in Ref. \cite{mall2017epsilon}, hyperbolic tangent function in Ref. \cite{taheri2018generic}, and normalized $L^2$-norm function; the initial guess for the shooting variable is randomly chosen in a defined domain $\boldsymbol{\eta}(0) = [\lambda_1(0), \lambda_2(0), t_f]^T$ as specified in Table~\ref{Table:InitialGuess}.
\begin{table}[!htp]
\centering
\caption{Initial guess for the shooting variable $\boldsymbol{\eta}(0)$}
\begin{tabular}{cccc}
\hline
Parameter  &$\lambda_1(0)$ &$\lambda_2(0)$  &$t_f$ \\
\hline
Initialization domain    &$[0,1]$    &$[0,1]$  &$[1,3]$ \\
\hline
\label{Table:InitialGuess}
\end{tabular}
\end{table}

The smoothing parameters for the epsilon-trig regularization, hyperbolic tangent function, and normalized $L^2$-norm function are set as $\epsilon = 1 \times 10^{-4}$, $\rho = 1\times 10^{-6}$, and $\delta = 1\times 10^{-8}$, respectively. The absolute and relative tolerances for propagating the state-costate equations 
are set to $1.0 \times 10^{-9}$. The default tolerance for \it{fsolve} \rm{ is set to $1.0 \times 10^{-9}$}.

\begin{figure}[!htp]
\centering
\includegraphics[width = 0.75\linewidth]{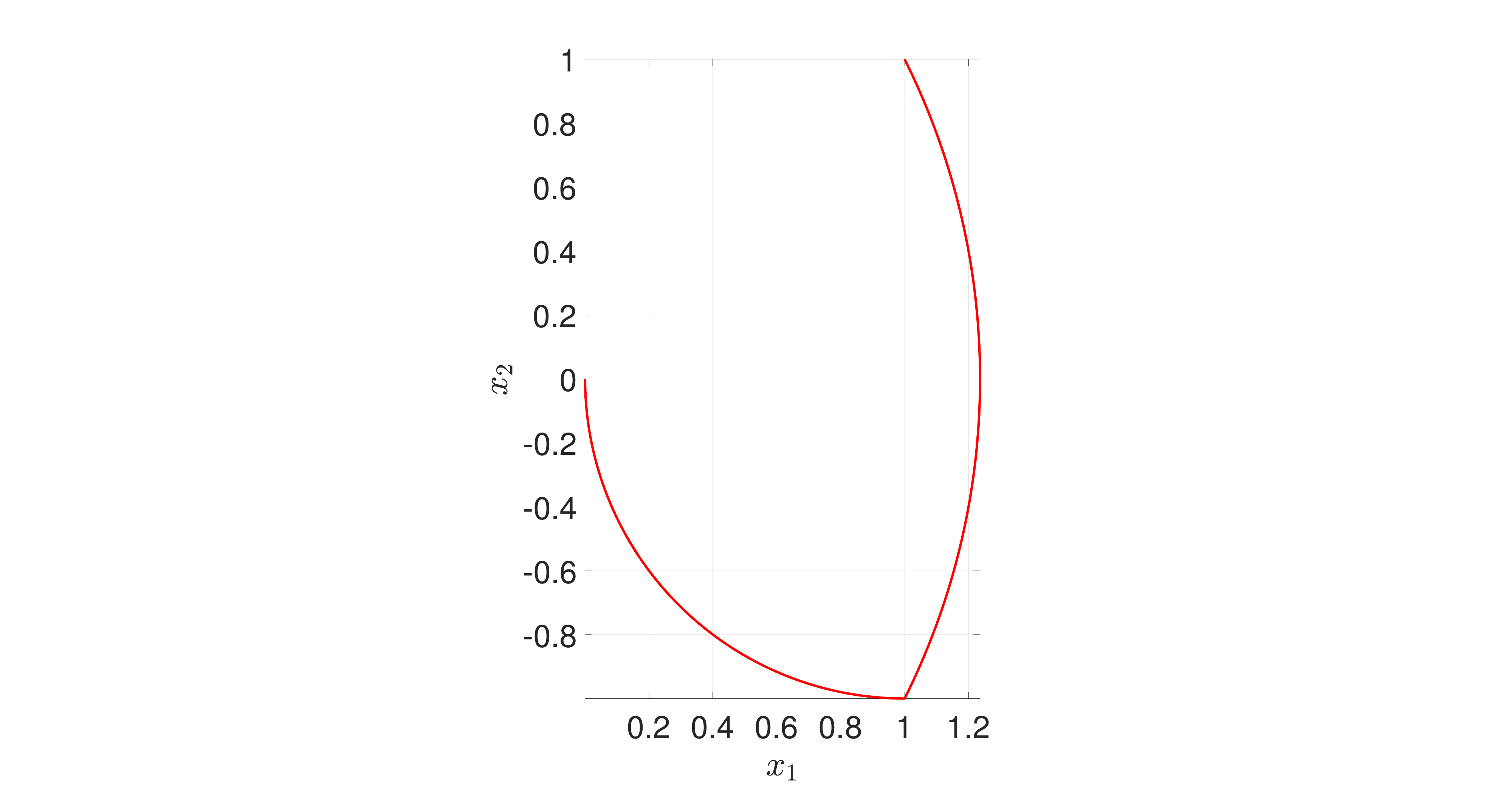}
\caption{Phase portrait for the minimal-time oscillator problem.}
\label{Fig:phase}
\end{figure}

\begin{figure}[!htp]
\centering
\includegraphics[width = 0.75\linewidth]{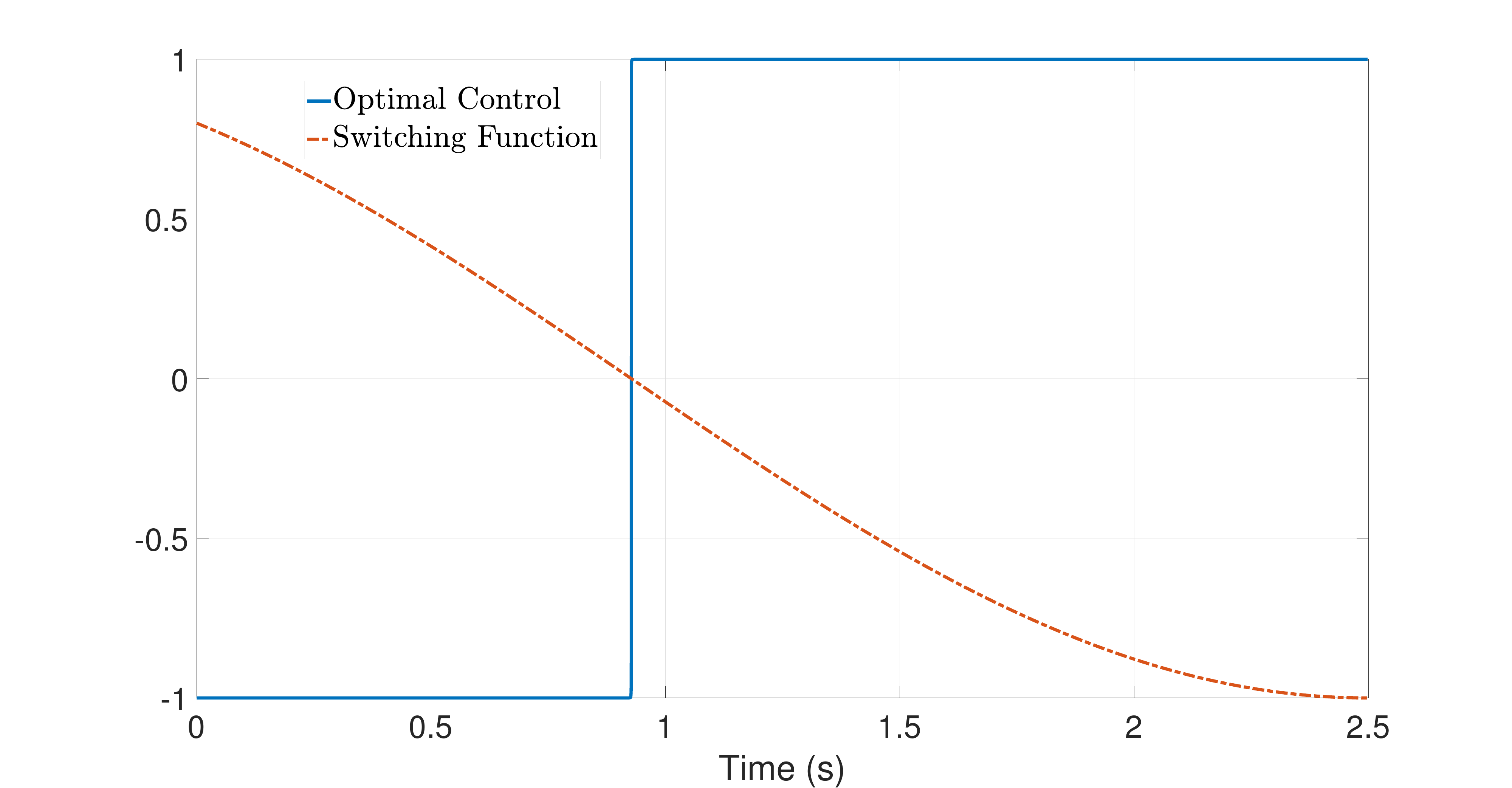}
\caption{Optimal control and switching function profiles.}
\label{Fig:control}
\end{figure}

The phase portrait, via using the normalized $L^2$-norm function as the smoothing method, is displayed in Fig.~\ref{Fig:phase}. 
The corresponding optimal control and switching function profiles are presented in Fig.~\ref{Fig:control}, and it foresees one switch for the optimal control at $t = 0.9273$ s. 
Table~\ref{Table:Results_1} summarizes the numerical results in terms of the number of converged cases, mean value of the cost function $t_f$, mean value of $\lvert\lvert\boldsymbol{\zeta}\rvert\rvert$ ($\lvert\lvert \cdot \rvert\rvert$ denotes the Euclidean norm), 
and computational time for three smoothing methods. 
The convergence rate for the normalized $L^2$-norm function $(95.58\%)$  is slightly lower than that of epsilon-trig regularization $(95.88\%)$, while the hyperbolic tangent function shows the lowest convergence rate, only $88.21\%$. The normalized $L^2$-norm function results in solutions with slightly less accuracies than another two smoothing methods, as shown by the mean values of resulting $t_f$ and $\lvert\lvert\boldsymbol{\zeta}\rvert\rvert$; however, the computational time for the normalized $L^2$-norm function is lower than those of other two smoothing methods. It is worth mentioning that the procedure for implementing the epsilon-trig regularization requires introducing an error control term, which would complicate the formulation of the TPBVP. On the other hand, both hyperbolic tangent function and normalized $L^2$-norm function can be readily implemented by approximating the bang-bang control using a smoothing function.
\begin{table}[!htp]
\centering
\caption{Comparison for three different smoothing methods} 
\begin{tabular}{ccccc}
\hline
Smoothing method  & No. of converged cases & Mean of $t_f$ (s)  & Mean of $\lvert\lvert\boldsymbol{\zeta}\rvert\rvert$ & Mean computational time (s)\\
\hline
Epsilon-trig regularization~\cite{mall2017epsilon}    &$9,588$    &$2.4980914$ &$4.8540 \times 10^{-9}$ &$0.12038$ \\
Hyperbolic tangent function~\cite{taheri2018generic}  &$8,821$    &$2.4980913$ &$7.0880 \times 10^{-9}$  &$0.04341$ \\
Normalized $L^2$-norm function       &$9,558$                  &$2.4980916$   &$8.4384 \times 10^{-8}$       &$0.03735$ \\
\hline
\label{Table:Results_1}
\end{tabular}
\end{table}
\subsection{Numerical Simulations on a Minimal-Fuel Low-Thrust Trajectory Optimization Problem}
\subsubsection{Problem Statement}
Here we consider an OCP for a spacecraft to transfer from a Geostationary Transfer Orbit (GTO) to the  Geosynchronous Equatorial Orbit  (GEO) while consuming minimal fuel. To avoid singularities from the classic orbit elements, the Modified Equinoctial Elements (MEEs) are employed to describe the spacecraft motion as \cite{walker1985set}
\begin{align}
\begin{cases}
\dot{\boldsymbol{x}}(t) = \boldsymbol{M}(\boldsymbol{x},t)\boldsymbol{\alpha}(t)\frac{u(t)T_{max}}{m(t)} + \boldsymbol{D}(\boldsymbol{x},t),\\
\dot{m}(t) = -\frac{u(t) T_{max}}{g_0 I_{sp}},
\label{EQ:dynamics_P2}
\end{cases}
\end{align}
in which $\boldsymbol{x} = [p,f,g,h,k,L]^T$ are the six equinoctial orbit
elements, $p$ the semilatus rectum, $[f,g]^T$ the eccentricity vector, $[h,k]^T$ the inclination vector, and $L$ is the true longitude; $m$ denotes the mass of the spacecraft; $T_{max}$ is the maximal thrust magnitude, $I_{sp}$ is the specific impulse of the engine, and $g_0$ is the gravitational acceleration at sea level. The engine throttle $u\in [0,1]$ and the unit vector of the thrust direction $\boldsymbol{\alpha}$ are the controls to be determined. The mathematical expressions for $\boldsymbol{M}(\boldsymbol{x})$ describing the effects of the thrust and $\boldsymbol{D}(\boldsymbol{x})$ denoting the control-free drift term, are presented in 
Refs.~\cite{gao2004low,gergaud2007orbital}.

The corresponding OCP aims to find the optimal controls $u^*(t)$ and $\boldsymbol{\alpha}^*(t)$ that drive the spacecraft from the GTO to GEO with final time $t_f$ fixed while minimizing the performance index 
\begin{align}
J = \frac{T_{max}}{g_0 I_{sp}}\int_{0}^{t_f} u(t)~ dt.
\end{align}

Before proceeding, it should be noted that normalization is usually required for better numerical conditioning; see, e.g., Ref. \cite{pan2019finding}. Let $\boldsymbol{\lambda} = [\lambda_p,\lambda_f,\lambda_g,\lambda_h,\lambda_k,\lambda_L]^T$ denote the costate vector related to the state vector $\boldsymbol{x}$, and let $\lambda_m$ represent the costate of mass. Then, the Hamiltonian can be defined as
\begin{align}
\mathscr{H}(t) = \frac{T_{max}}{g_0 I_{sp}} u(t) + \boldsymbol{\lambda}^T(t) [\boldsymbol{M}(\boldsymbol{x,t})\boldsymbol{\alpha}(t)\frac{u(t)T_{max}}{m(t)} + \boldsymbol{D}(\boldsymbol{x},t)] - \lambda_m(t) \frac{u(t) T_{max}}{g_0 I_{sp}}.
\label{EQ:Ham_P2}
\end{align}
The corresponding equations of the costates are expressed as \cite{Pontryagin}
\begin{align}
\begin{cases}
\dot{\boldsymbol{\lambda}}(t) = -[\frac{\partial \mathscr{H}(t)}{\partial \boldsymbol{x}(t)}]^T,\\
\dot{\lambda}_m(t) = -\frac{\partial \mathscr{H}(t)}{\partial m(t)}.
\label{EQ:costatedynamics}
\end{cases}
\end{align}
At the initial time $t_0 = 0$, the initial condition is given as $\boldsymbol{x}(0) = [p_0,f_0,g_0,h_0,k_0,L_0]^T$ ; at the fixed final time $t_f$, the final condition for the spacecraft is set as
\begin{equation}
\begin{aligned}
\boldsymbol{\Phi} (\boldsymbol{x} (t_f), t_f) =
[p(t_f) - p_f, f(t_f) -f_f,g(t_f) - g_f,h(t_f) -h_f,k(t_f) - k_f]^T = \boldsymbol{0},
\label{EQ:terminal}
\end{aligned}
\end{equation}
where $p_f$, $f_f$, $g_f$, $h_f$, $k_f$ are the specified MEEs of the target orbit.  Additionally, the true longitude $L_f$ and mass $m$ at $t_f$ are left free, thus we have
\begin{align}
\begin{cases}
\lambda_L(t_f) = 0, \\
\lambda_m(t_f) = 0.
\label{EQ:costatemass}
\end{cases}
\end{align}
Notice that the optimal engine throttle takes a bang-bang control form for a fuel-optimal low-thrust transfer problem \cite{haberkorn2004low}. Then, the optimal thrust direction $\boldsymbol{\alpha}^*(t)$ and engine throttle $u^*(t)$ should satisfy  
\begin{align}
\boldsymbol{\alpha}^*(t) = -\frac{\boldsymbol{M^T}(\boldsymbol{x},t)\boldsymbol{\lambda}(t)}
{\lvert \lvert \boldsymbol{M^T}(\boldsymbol{x},t)\boldsymbol{\lambda(t)} \rvert \rvert} 
\label{EQ:optimalControl_alpha}
\end{align}
and
\begin{align}
u^*(t) = 
\begin{cases}
1, \rm{if}~\it{S(t) < \rm{0}}, \\
0, \rm{if}~\it{S(t) > \rm{0}},
\end{cases}
\label{EQ:optimalControl_thrust}
\end{align}
where $S(t)$ denotes the switching function defined as
\begin{align}
S(t) = 1 - \frac{g_0 I_{sp} {\lvert \lvert \boldsymbol{M^T}(\boldsymbol{x},t)\boldsymbol{\lambda}(t) \rvert \rvert}}{m(t)} - \lambda_m(t).
\label{EQ:optimalControl_S}
\end{align}

We use the normalized $L^2$-norm function to approximate the discontinuous optimal control in Eq.~(\ref{EQ:optimalControl_thrust}), i.e.,
\begin{align}
u^*(t) \approx u^*(t,\delta) = \frac{1}{2}(1-\frac{S(t)}{\sqrt{\delta+\lvert S(t)\rvert^2}}).
\label{EQ:smoothing_P2}
\end{align}
Thus, addressing the responding OCP is equivalent to solving the shooting function
\begin{align}
\boldsymbol{\zeta}(\boldsymbol{\eta}(0)) = [\boldsymbol{\Phi} (\boldsymbol{x} (t_f), t_f), \lambda_L(t_f), \lambda_m(t_f)] = \boldsymbol{0},
\label{EQ:shooting_function_p2}
\end{align}
where 
$\boldsymbol{\eta}(0) = [\lambda_p(0), \lambda_f(0), \lambda_g(0),\lambda_h(0),\lambda_k(0),\lambda_L(0),\lambda_m(0)]^T$ is the shooting variable.
\subsubsection{Numerical Results}
Consider a spacecraft with an initial mass of $1,500$ kg, and specific impulse $I_{sp}$ of $2,000$ s. The maximal thrust magnitude $T_{max}$ is $1.0$ N, and the constant $g_0$ is equal to $9.80665$ $ \rm{m/s^2}$. 
In this case, the initial thrust-to-weight ratio of the spacecraft is $6.7981 \times 10^{-5}$, which can be characterized as a very low-thrust trajectory optimization problem, making it quite challenging to solve \cite{pan2019finding}. For the initial elliptical orbit and final circular orbit, we use the MEEs in Refs. \cite{zhang2023solution,pan2019finding} as presented in Table~\ref{Table:MEEs}. The final transfer time is fixed as $1,000$ hours \cite{pan2019finding}.
\begin{table}[!htp]
\centering
\caption{MEEs of the GTO and GEO}
\begin{tabular}{ccccccc}
\hline
MEEs  &$p$~($km$) &$f$  &$g$  &$h$ &$k$  &$L$~(rad)\\ 
\hline
GTO  & $11,623$    & $0.75$  & $0$ & $0.0612$ & $0$ & $\pi$\\ 
GEO & $42,165$    & $0$  & $0$ & $0$ & $0$ & free\\ 
\hline
\label{Table:MEEs}
\end{tabular}
\end{table}

Analogously, $100$ simulations are performed for the smoothing technique in Ref. \cite{taheri2018generic} and the proposed method. The initial guess for the shooting variable in both methods is randomly generated in a defined domain $\boldsymbol{\eta}(0) = [\lambda_p(0), \lambda_f(0), \lambda_g(0),\lambda_h(0),\lambda_k(0),\lambda_L(0),\lambda_m(0)]^T$ as shown in Table~\ref{Table:InitialGuess_problem2}.
\begin{table}[!htp]
\centering
\caption{Initial guess for the shooting variable $\boldsymbol{\eta}(0)$}
\begin{tabular}{cccccccc}
\hline
Parameter  &$\lambda_p(0)$ &$\lambda_f(0)$  &$\lambda_g(0)$ 
&$\lambda_h(0)$ &$\lambda_k(0)$  &$\lambda_L(0)$ &$\lambda_m(0)$\\
\hline
Initialization domain    &$[0,0.1]$    &$[0,0.1]$  &$[0,0.1]$ &$[0,0.1]$    &$[0,0.1]$  &$[0,0.1]$ &$[0,0.1]$\\
\hline
\label{Table:InitialGuess_problem2}
\end{tabular}
\end{table}

To facilitate the convergence of the indirect method, a continuation procedure is devised as follows. Recall that $\rho$ denotes the smoothing constant for the hyperbolic tangent function in Ref. \cite{taheri2018generic}, and $\delta$ represents the smoothing constant for the technique proposed in this paper.
\begin{algorithm}[htb]
\setstretch{1.35} 
\caption{Algorithm for the continuation procedure}
\label{alg:Framwork}
(1) Set the absolute and relative tolerances for propagating the state-costate dynamics as $1.0 \times 10^{-12}$. Set the default tolerance for \it{fsolve} \rm{as $1.0 \times 10^{-6}$}.\\
(2) Let $\rho_0 = 1~(\delta_0 = 1)$. Set $\rho_{min} = 1 \times 10^{-6}~(\delta_{min} = 1 \times 10^{-8})$. Let $n = 0$. \\
(3) Initialize $\boldsymbol{\eta}(0)$ by Table~\ref{Table:InitialGuess_problem2}. \\
(4) Solve the TPBVP corresponding to $\rho_{n} ~(\delta_{n})$. \\
(5) If the indirect method converges, let $n = n + 1$ and go to step 6. If not, perturbate $\boldsymbol{\eta}(0)$ by $\boldsymbol{\eta}(0) = \boldsymbol{\eta}(0) + \frac{\rho}{100} ~(\frac{\delta}{100})~\times$ rand(7,1), go to step 4. If the indirect method does not converge after $10$ perturbations over $\boldsymbol{\eta}(0)$, the indirect method fails and exit. \\
(6) Let $\rho_{n} = \frac{\rho_{n}}{10}  ~(\delta_{n}= \frac{\delta_{n}}{10})$. Initialize $\boldsymbol{\eta}(0)$ using the preceding convergent solution that corresponds to $\rho_{n} ~(\delta_{n})$; go to step 4. \\
(7) Output the convergent shooting variable as well as the optimal solution once $\rho_{n} \leq \rho_{min}$~$(\delta_{n} \leq \delta_{min})$.
\end{algorithm}

Notice that for a very low-thrust geocentric transfer problem, numerous local optimal solutions may exist \cite{pan2019finding}. Table~\ref{Table:Results_2}
\begin{table}[!htp]
\centering
\caption{Comparison for two different smoothing methods} 
\begin{tabular}{ccc}
\hline
Smoothing method  & Hyperbolic tangent function~\cite{taheri2018generic}
& Normalized $L^2$-norm function \\
\hline
No. of converged cases &$36$    &$42$ \\
$\Delta m$ (kg)  &$[133.6921, 183.5489]$  &$[133.6844,  183.5489]$\\
No. of revolutions   &$[48,75]$  &$[47, 75]$ \\  
Mean computational time (s)  &$187.9348$  &$208.6674$\\
Mean of  $\lvert\lvert\boldsymbol{\zeta}\rvert\rvert$ &$2.1157\times 10^{-6}$  &$2.4168 \times 10^{-6}$\\
\hline
\label{Table:Results_2}
\end{tabular}
\end{table} 
demonstrates the results for two smoothing methods in terms of the number of converged cases, range of fuel consumption $\Delta m$, range of number of revolutions, mean computational time and accuracy $\lvert\lvert\boldsymbol{\zeta}\rvert\rvert$ in satisfying the shooting function Eq.~(\ref{EQ:shooting_function_p2}). It is worth mentioning that due to numerical errors, it is difficult to determine the number of locally optimal solutions. Based on the results for $\Delta m$ and number of revolutions, it is reasonable to conclude that both smoothing techniques find many locally optimal solutions. More importantly, the convergence rate for the normalized $L^2$-norm  function $(42\%)$ is slightly higher than that of the hyperbolic tangent function $(36\%)$, although it results in solutions with slightly lower accuracies in satisfying the terminal constraints. Additionally, due to the fact the smoothing constant $\delta$ for the normalized $L^2$-norm  function is smaller than the smoothing constant $\rho$ for the hyperbolic tangent function, the normalized $L^2$-norm  function needs more steps in the continuation procedure, resulting in higher computational time. 

\begin{figure}[!htp]
\centering
\includegraphics[width = 1\linewidth]{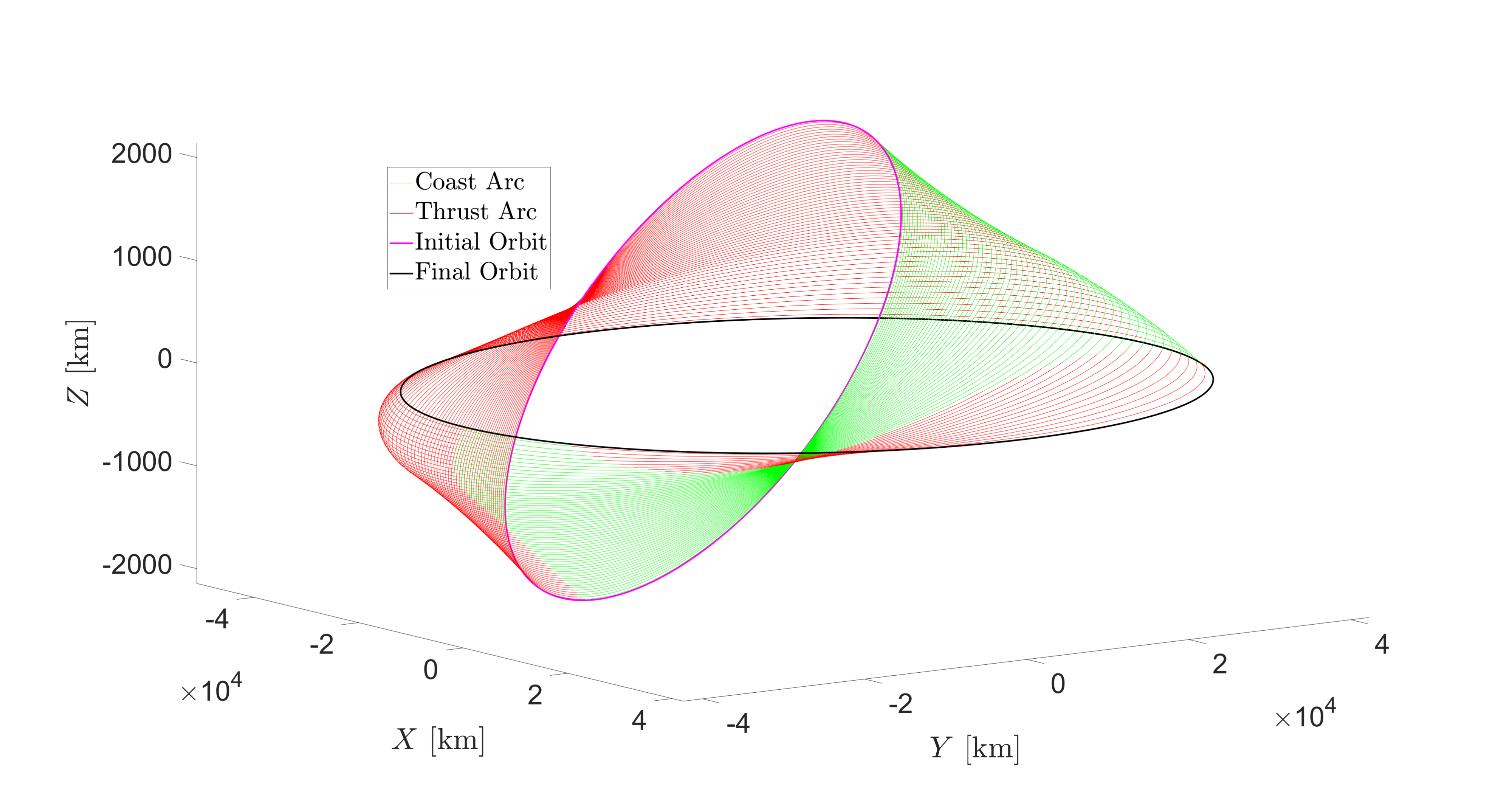}
\caption{Three-dimensional transfer trajectory of the potentially best solution (The axis $Z$ is not to scale).}
\label{Fig:low_thrust_3d}
\end{figure}

\begin{figure}[!htp]
\centering
\includegraphics[width = 0.9\linewidth]{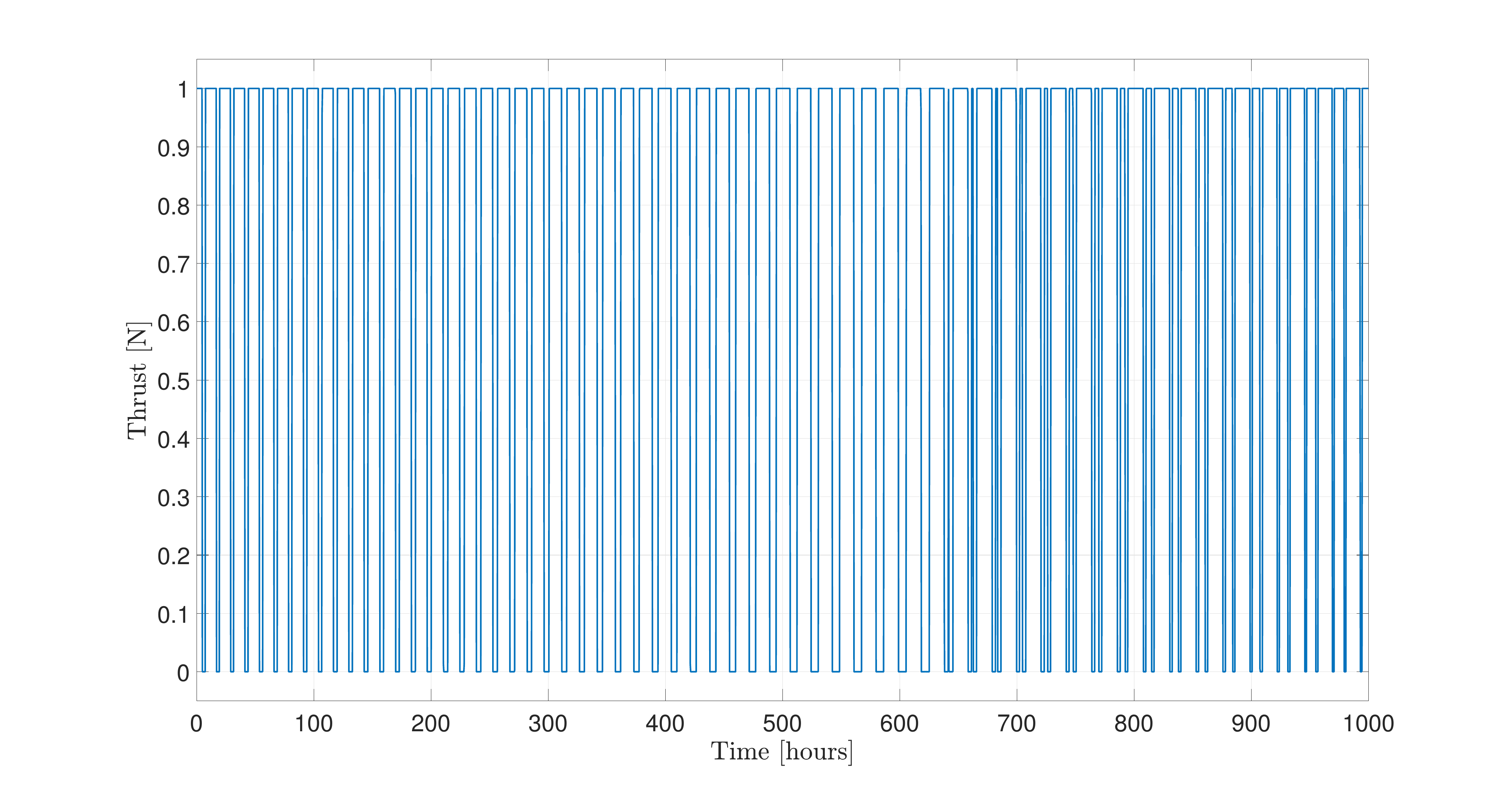}
\caption{Engine thrust profile of the potentially best solution.}
\label{Fig:engine_throttle}
\end{figure}

\begin{figure}[!htp]
\centering
\includegraphics[width = 1\linewidth]{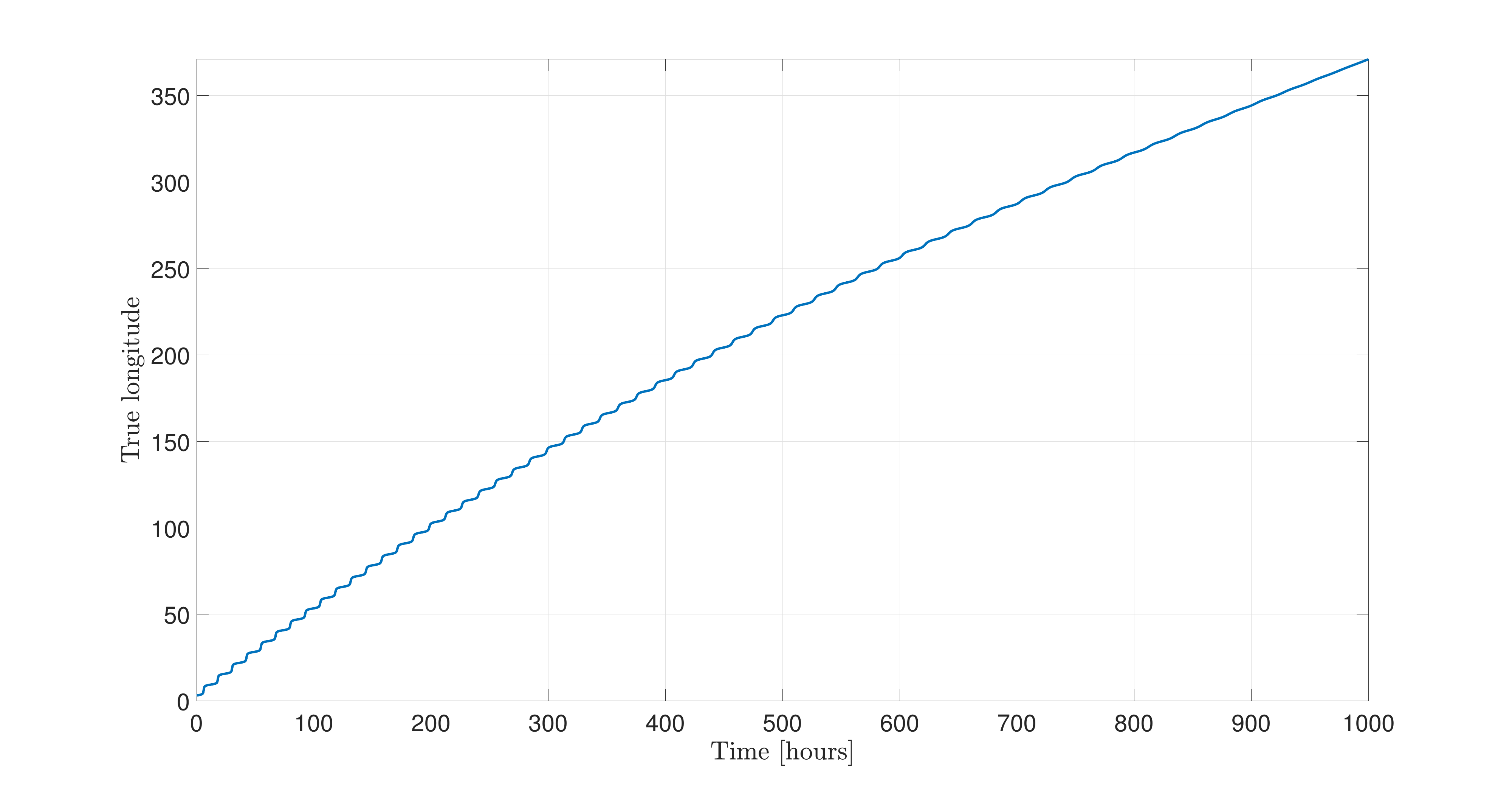}
\caption{True longitude profile of the potentially best solution.}
\label{Fig:true longitude}
\end{figure}

Among the obtained solutions, the potentially best solution can be identified with the overall minimal fuel consumption of $133.6844$ kg as shown in Table~\ref{Table:Results_2}, and this result has been previously demonstrated in Ref. \cite{pan2019finding}. 
\begin{figure}[!htp]
\centering
\includegraphics[width = 1\linewidth]{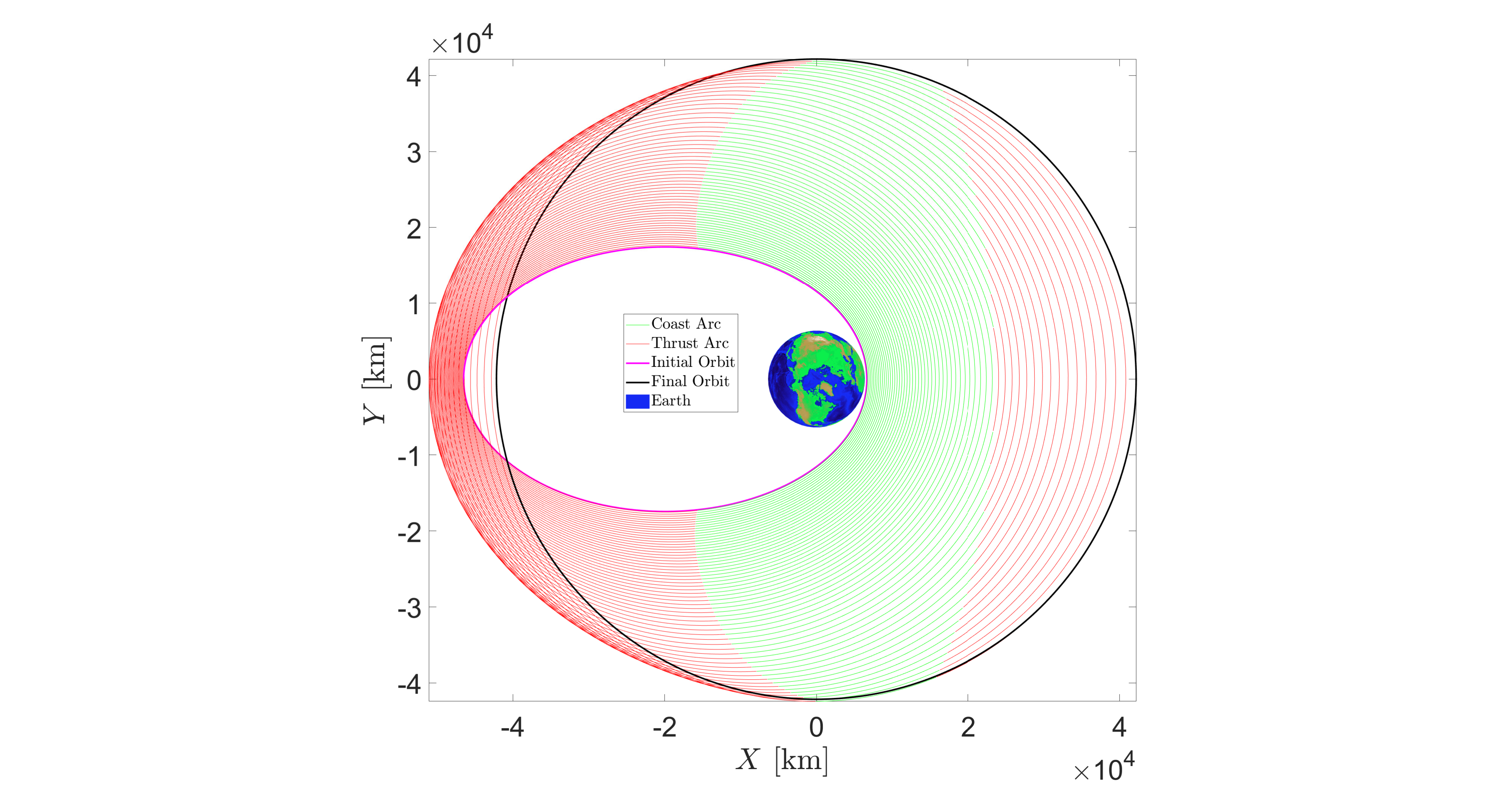}
\caption{Transfer trajectory of the potentially best solution projected onto the $XY$ plane.}
\label{Fig:low_thrust}
\end{figure}
Fig.~\ref{Fig:low_thrust} shows the minimal-fuel low-thrust transfer trajectory projected onto the $XY$ plane, and the three-dimensional transfer trajectory from GTO to GEO is illustrated in Fig.~\ref{Fig:low_thrust_3d} (The Cartesian coordinate system $XYZ$ follows the right-hand rule. The plane $XY$ and the circular target orbit are coplanar, and its origin is coincident with the center of the circular target orbit). 
The engine thrust profile during the transfer is displayed in Fig.~\ref{Fig:engine_throttle}, from which it can be seen that the engine throttle involves a total of $149$ switches. 
Fig.~\ref{Fig:true longitude} depicts the true longitude profile during the transfer with a total of $INT(\frac{L_f - L_0}{2\pi}) = 58$ revolutions ($INT$ means rounding down to the nearest integer), which has also been  reported in Ref. \cite{pan2019finding}. The results show that the proposed smoothing method can be readily applied to finding the potentially best solution among numerous locally optimal solutions for a minimal-fuel low-thrust trajectory optimization problem that involves many revolutions.
\section{Conclusions} \label{Colus}
In this work, a bounded smooth function, called normalized $L^2$-norm function is proposed to smooth bang-bang controls. 
Unlike some smoothing or regularization techniques that modify the cost function or state equation, the proposed smooth function approximates the sign function in terms of the switching function, leading to a smooth optimal control profile. As a consequence, the proposed method can be readily embedded into the indirect method once the switching function is derived from the Pontryagin’s Minimum Principle. The application of the proposed method is demonstrated by two different problems: 1) a minimal-time oscillator problem, and 2) a minimal-fuel low-thrust trajectory optimization problem that involves many revolutions. 

\section*{Acknowledgement}
This research was supported by the Key Research and Development Program of Zhejiang Province (No. 2020C05001).

\bibliography{AIAA}

\end{document}